\documentclass[12pt]{article}
\usepackage{amssymb,amsmath}
\usepackage{cases}
\usepackage{amsfonts}
\usepackage{color,xcolor,cite}
\usepackage[left=2.0cm,right=2.0cm,top=2.0cm,bottom=2.0cm]{geometry}
\usepackage[colorlinks,citecolor=blue,urlcolor=blue]{hyperref}

\newtheorem{theorem}{Theorem}[section]

\newtheorem{lemma}{Lemma}[section]
\newtheorem{proposition}{Proposition}[section]

\newtheorem{definition}{Definition}[section]
\newtheorem{remark}{Remark}[section]

\newcommand{\bal}{\begin{align}}
\newcommand{\bbal}{\begin{align*}}
\newcommand{\beq}{\begin{equation}}
\newcommand{\eeq}{\end{equation}}
\newcommand{\bca}{\begin{cases}}
\newcommand{\eca}{\end{cases}}

\newcommand{\pa}{\partial}
\newcommand{\fr}{\frac}

\newcommand{\De}{\Delta}

\newcommand{\ep}{\varepsilon}
\newcommand{\dd}{\mathrm{d}}

\newcommand{\R}{\mathbb{R}}
\newcommand{\vv}{\mathbf{v}}

\newcommand{\bi}{\Big}

\begin{document}
\title{Ill-posedness for the Camassa-Holm and related equations in Besov spaces }

\author{Jinlu Li$^{1}$, Yanghai Yu$^{2,}$\footnote{E-mail: lijinlu@gnnu.edu.cn; yuyanghai214@sina.com(Corresponding author); mathzwp2010@163.com} and Weipeng Zhu$^{3}$\\
\small $^1$ School of Mathematics and Computer Sciences, Gannan Normal University, Ganzhou 341000, China\\
\small $^2$ School of Mathematics and Statistics, Anhui Normal University, Wuhu 241002, China\\
\small $^3$ School of Mathematics and Big Data, Foshan University, Foshan, Guangdong 528000, China}

\date{\today}

\maketitle\noindent{\hrulefill}

{\bf Abstract:} In this paper, we give a construction of $u_0\in B^\sigma_{p,\infty}$ such that
corresponding solution to the Camassa-Holm equation starting from $u_0$ is
discontinuous at $t = 0$ in the metric of $B^\sigma_{p,\infty}$, which implies the ill-posedness for this equation in $B^\sigma_{p,\infty}$. We also apply our method to the b-equation and Novikov equation.

{\bf Keywords:} Camassa-Holm equation; Shallow water wave models; Ill-posedness; Besov space.

{\bf MSC (2010):} 35Q53, 37K10.
\vskip0mm\noindent{\hrulefill}

\section{Introduction}\label{sec1}
In this paper, we consider the question of the well-posedness of the Cauchy problems
to a class of shallow water wave equations on the real-line, such as the Novikov equation and the b-equation contains two integrable members, the Camassa-Holm equation and the Degasperis-Procesi equation. The method we used in this paper is very simple and can be applied equally well to these equations. Thus, in order to elucidate the main ideas, our attention in this paper will be focused on the Camassa-Holm equation. For other equations, we give the results as remarks.

The Cauchy problem for the well-known Camassa-Holm equation reads as follows
\begin{align}\label{ch0}
\begin{cases}
u_t-u_{xxt}+3uu_x=2u_xu_{xx}+uu_{xxx}, \; &(x,t)\in \R\times\R^+,\\
u(x,t=0)=u_0,\; &x\in \R.
\end{cases}
\end{align}
Here the scalar function $u = u(t, x)$ stands for the fluid velocity at time $t\geq0$ in the $x$ direction.
\eqref{ch0} was firstly proposed in the context of hereditary symmetries studied in \cite{Fokas} and then was derived explicitly as a water wave equation by Camassa--Holm \cite{Camassa}. \eqref{ch0} is completely integrable \cite{Camassa,Constantin-P} with a bi-Hamiltonian structure \cite{Constantin-E,Fokas} and infinitely many conservation laws \cite{Camassa,Fokas}. Also, it admits exact peaked
soliton solutions (peakons) of the form $ce^{-|x-ct|}$ with $c>0$, which are orbitally stable \cite{Constantin.Strauss} and models wave breaking (i.e., the solution remains bounded, while its slope becomes unbounded in finite time \cite{Constantin,Escher2,Escher3}). The peaked solitons present the characteristic for the travelling water waves of greatest height and largest amplitude and arise as solutions to the free-boundary problem for incompressible Euler equations over a flat bed, see Refs. \cite{Constantin-I,Escher4,Escher5,Toland} for the details.

Because of the interesting and remarkable features as mentioned above, the Camassa-Holm equation has attracted much attention as a class of integrable shallow water wave equations in recent twenty years. Its systematic mathematical study was initiated in a series of papers by Constantin and Escher, see \cite{Escher1,Escher2,Escher3,Escher4,Escher5}.
There is an extensive literature about the well-posedness theory of the Camassa-Holm equation. Before stating related results for \eqref{ch0} precisely, we recall the notion of well-posedness in the sense of Hadamard. We say that the Cauchy problem \eqref{ch0} is Hadamard (locally) well-posed in a Banach space $X$ if for any data $u_0\in X$ there exists (at least for a short time) $T>0$ and a unique solution in the space $\mathcal{C}([0,T),X)$ which depends continuously on the data. In particular, we say that the solution map is
continuous if for any $u_0\in X$, there exists a neighborhood $B \subset X$ of $u_0$ such that
for every $u \in B$ the map $u \mapsto U$ from $B$ to $\mathcal{C}([0, T]; X)$ is continuous, where $U$
denotes the solution to \eqref{ch0} with initial data $u_0$.

Next we only review some results concerning the well-posedness of the Cauchy problem \eqref{ch0} (see
also the survey \cite{Molinet}). Li--Olver \cite{li2000} proved that the Cauchy problem \eqref{ch0} is
locally well-posed with the initial data $u_0\in H^s(\R)$ with $s > 3/2$ (see also \cite{GB}). Danchin \cite{d1} proved
the local existence and uniqueness of strong solutions to \eqref{ch0} with initial data in $B^s_{p,r}$ for $s > \max\{1 + 1/p , 3/2\}$ with $p\in[1,\infty]$ and $r\in[1,\infty)$. Meanwhile, he \cite{d1} only obtained the continuity of the solution map of \eqref{ch0} with respect to the initial data in the space $\mathcal{C}([0, T ];B^{s'}_{p,r})$ with any $s'<s$. Li--Yin \cite{Li-Yin1} proved the continuity of the solution map of \eqref{ch0} with respect to the initial data in the space $\mathcal{C}([0, T];B^{s}_{p,r})$ with $r<\infty$. In particular, they \cite{Li-Yin1} proved that the solution map of \eqref{ch0} is weak continuous with respect to initial data $u_0\in B^s_{p,\infty}$. For the endpoints, Danchin \cite{d3} obtained the local well-posedness in the space $B^{3/2}_{2,1}$ and ill-posedness in $B^{3/2}_{2,\infty}$ (the data-to-solution map is not continuous by using peakon solution). In our recent papers\cite{Li1, Li2}, we proved the non-uniform dependence on initial data for \eqref{ch0} under both the framework of Besov spaces $B^s_{p,r}$ for $s>\max\big\{1+1/p, 3/2\big\}$ with $p\in[1,\infty], r\in[1,\infty)$ and $B^{3/2}_{2,1}$ (see \cite{H-M,H-K,H-K-M} for earlier results in $H^s$). Byers \cite{Byers} proved that the Camassa-Holm equation is ill-posed in $H^s$ for $s < 3/2$ in the sense of norm inflation, which means that $H^{3/2}$ is the critical Sobolev space for well-posedness. Moreover, for the intermediate cases, Guo et al.\cite{Guo-Yin} established the ill-posedness of \eqref{ch0} in the critical Sobolev
space $H^{3/2}$ and even in the Besov space $B_{p,r}^{1+1/p}$ with $p\in[1,\infty],r\in(1,\infty]$ by proving the norm inflation.

{\emph{Question appears: Whether or not the continuity of the data-to-solution map with values in $L^\infty(0, T;B^s_{p,\infty})$ with $s>\max\{1+1/p, 3/2\}$ and $1\leq p\leq\infty$ holds for the Camassa-Holm equation. To the best of our knowledge, this issue has not been solved yet, which is our goal in this paper.}}

Before stating our main result, we transform \eqref{ch0} equivalently into the following nonlinear transport type equation
\begin{align}\label{ch}
\begin{cases}
\partial_tu+u\pa_xu=-\pa_x(1-\pa^2_x)^{-1}\big(u^2
+\frac{1}{2}u^2_x\big):=\mathbf{P}(u), \; &(x,t)\in \R\times\R^+,\\
u(x,t=0)=u_0,\; &x\in \R.
\end{cases}
\end{align}

Now let us state our main result of this paper.
\begin{theorem}\label{th1}
Let $\sigma>2+\max\big\{3/2,1+1/p\big\}$ with $1\leq p\leq \infty$. There exits $u_0\in B^\sigma_{p,\infty}(\R)$ and a positive constant $\ep_0$ such that the data-to-solution map $u_0\mapsto \mathbf{S}_{t}(u_0)$ of the Cauchy problem \eqref{ch} satisfies\bbal
\limsup_{t\to0^+}\|\mathbf{S}_{t}(u_0)-u_0\|_{B^\sigma_{p,\infty}}\geq \ep_0.
\end{align*}
\end{theorem}
\begin{remark}
Theorem \ref{th1} demonstrates the ill-posedness of the Camassa-Holm equation in $B^\sigma_{p,\infty}$. More precisely, there exists a $u_0\in B^\sigma_{p,\infty}$ such that corresponding solution to the CH equation that starts from $u_0$ does not converge back to $u_0$ in the metric of $B^\sigma_{p,\infty}$ as time goes to zero. Our key argument is to present a construction of initial data $u_0$.
\end{remark}
\begin{remark}
Theorem \ref{th1} also holds for the following b-family equation (see \cite{03A,03B,DP,16H} etc.):
\begin{align}\label{b-family}
\begin{cases}
\pa_tu+uu_x=-\pa_x(1-\pa^2_x)^{-1}\big(\frac{b}{2}u^2
+\frac{3-b}{2}u^2_x\big), &\quad (t,x)\in \R^+\times\R,\\
u(0,x)=u_0(x), &\quad x\in \R,
\end{cases}
\end{align}
\end{remark}
It should be mentioned that the Camassa-Holm equation corresponds to $b = 2$ and Degasperis-Procesi equation corresponds to $b = 3$.

The Cauchy problem for the Novikov equation reads as (see \cite{H-H,HHK2018,Li2,Li3} etc.)
\begin{align}\label{nov}
\begin{cases}
u_t+u^2u_x=-(1-\pa^2_x)^{-1}\Big(\frac12u_x^3
+\pa_x\big(\frac32uu^2_x+u^3\big)\Big)=:\mathbf{Q}(u), \\
u(0,x)=u_0,
\end{cases}
\end{align}
Then, we have the following
\begin{theorem}\label{th2}
Let $\sigma>\frac72$. There exits $u_0\in B^\sigma_{2,\infty}(\R)$ and a positive constant $\ep_0$ such that the data-to-solution map $u_0\mapsto \mathbf{S}_{t}(u_0)$ of \eqref{nov}
satisfies
\bbal
\limsup_{t\to0^+}\|\mathbf{S}_{t}(u_0)-u_0\|_{B^\sigma_{2,\infty}}\geq \ep_0.
\end{align*}
\end{theorem}

\begin{remark}
The difficulty in proving that Theorem \ref{th2} holds in $B^\sigma_{p,\infty}$ with $p\neq2$ consists mainly in the construction of initial data $u_0$ due to the appearance of $u^2u_x$, we do not pursue the issue in the current paper. However, for $p=2$, the initial data $u_0$ as constructed in Lemma \ref{le6} seems more simple.
\end{remark}

The remaining of this paper is organized as follows. In Section \ref{sec2} we list some notations and recall known results. In Section \ref{sec3} we present the proof of Theorem \ref{th1} by establishing some technical lemmas and propositions. In Section \ref{sec4} we prove Theorem \ref{th2}.

\section{Preliminaries}\label{sec2}
\subsection{General Notation}\label{sec2.1}
In the following, we denote by $\star$ the convolution.
Given a Banach space $X$, we denote its norm by $\|\cdot\|_{X}$. For $I\subset\R$, we denote by $\mathcal{C}(I;X)$ the set of continuous functions on $I$ with values in $X$. Sometimes we will denote $L^p(0,T;X)$ by $L_T^pX$.
For all $f\in \mathcal{S}'$, the Fourier transform $\mathcal{F}f$ (also denoted by $\widehat{f}$) is defined by
$$
\mathcal{F}f(\xi)=\widehat{f}(\xi)=\int_{\R}e^{-ix\xi}f(x)\dd x \quad\text{for any}\; \xi\in\R.
$$
\subsection{Littlewood-Paley Analysis}\label{sec2.2}
Next, we will recall some facts about the Littlewood-Paley decomposition, the nonhomogeneous Besov spaces and their some useful properties.
\begin{proposition}[Littlewood-Paley decomposition, See \cite{B.C.D}] Let $\mathcal{B}:=\{\xi\in\mathbb{R}:|\xi|\leq \frac 4 3\}$ and $\mathcal{C}:=\{\xi\in\mathbb{R}:\frac 3 4\leq|\xi|\leq \frac 8 3\}.$
There exist two radial functions $\chi\in C_c^{\infty}(\mathcal{B})$ and $\varphi\in C_c^{\infty}(\mathcal{C})$ both taking values in $[0,1]$ such that
\begin{align*}
&\chi(\xi)+\sum_{j\geq0}\varphi(2^{-j}\xi)=1 \quad \forall \;  \xi\in \R^d,\\
&\frac{1}{2} \leq \chi^{2}(\xi)+\sum_{j \geq 0} \varphi^{2}(2^{-j} \xi) \leq 1\quad \forall \;  \xi\in \R^d.
\end{align*}
\end{proposition}
\begin{remark}
It is easy to show that $\varphi(\xi)\equiv 1$ for $\frac43\leq |\xi|\leq \frac32$.
\end{remark}

\begin{definition}[See \cite{B.C.D}]
For every $u\in \mathcal{S'}(\mathbb{R})$, the Littlewood-Paley dyadic blocks ${\Delta}_j$ are defined as follows
\begin{numcases}{\Delta_ju=}
0, & if $j\leq-2$;\nonumber\\
\chi(D)u=\mathcal{F}^{-1}(\chi \mathcal{F}u), & if $j=-1$;\nonumber\\
\varphi(2^{-j}D)u=\mathcal{F}^{-1}\big(\varphi(2^{-j}\cdot)\mathcal{F}u\big), & if $j\geq0$.\nonumber
\end{numcases}
In the inhomogeneous case, the following Littlewood-Paley decomposition makes sense
$$
u=\sum_{j\geq-1}{\Delta}_ju,\quad \forall\;u\in \mathcal{S'}(\mathbb{R}).
$$
\end{definition}
\begin{definition}[See \cite{B.C.D}]
Let $s\in\mathbb{R}$ and $(p,r)\in[1, \infty]^2$. The nonhomogeneous Besov space $B^{s}_{p,r}(\R)$ is defined by
\begin{align*}
B^{s}_{p,r}(\R):=\Big\{f\in \mathcal{S}'(\R):\;\|f\|_{B^{s}_{p,r}(\mathbb{R})}<\infty\Big\},
\end{align*}
where
\begin{numcases}{\|f\|_{B^{s}_{p,r}(\mathbb{R})}=}
\left(\sum_{j\geq-1}2^{sjr}\|\Delta_jf\|^r_{L^p(\mathbb{R})}\right)^{\fr1r}, &if $1\leq r<\infty$,\nonumber\\
\sup_{j\geq-1}2^{sj}\|\Delta_jf\|_{L^p(\mathbb{R})}, &if $r=\infty$.\nonumber
\end{numcases}
\end{definition}
\begin{remark}\label{re3}
It should be emphasized that the following embedding will be often used implicity:
$$B^s_{p,q}(\R)\hookrightarrow B^t_{p,r}(\R)\quad\text{for}\;s>t\quad\text{or}\quad s=t,1\leq q\leq r\leq\infty.$$
\end{remark}
Finally, we give some important properties which will be also often used throughout the paper.
\begin{lemma}[See \cite{B.C.D}]\label{le1}
Let $(p,r)\in[1, \infty]^2$ and $s>\max\big\{1+\frac1p,\frac32\big\}$. Then we have
\bbal
&\|uv\|_{B^{s-2}_{p,r}(\R)}\leq C\|u\|_{B^{s-2}_{p,r}(\R)}\|v\|_{B^{s-1}_{p,r}(\R)}.
\end{align*}
Hence, for the terms $\mathbf{P}(u)$ and $\mathbf{P}(v)$, we have
\bbal
&\|\mathbf{P}(u)-\mathbf{P}(v)\|_{B^{s-1}_{p,r}(\R)}\leq C\|u-v\|_{B^{s-1}_{p,r}(\R)}\|u+v\|_{B^{s}_{p,r}(\R)}.
\end{align*}
\end{lemma}
\begin{lemma}[See \cite{B.C.D}]\label{le2}
For $(p,r)\in[1, \infty]^2$ and $s>0$, $B^s_{p,r}(\R)\cap L^\infty(\R)$ is an algebra. Moreover, $B^{\frac{1}{p}}_{p,1}(\R)\hookrightarrow L^\infty(\R)$, and for any $u,v \in B^s_{p,r}(\R)\cap L^\infty(\R)$, we have
\bbal
&\|uv\|_{B^{s}_{p,r}(\R)}\leq C\big(\|u\|_{B^{s}_{p,r}(\R)}\|v\|_{L^\infty(\R)}+\|v\|_{B^{s}_{p,r}(\R)}\|u\|_{L^\infty(\R)}\big).
\end{align*}
\end{lemma}

\section{Proof of Theorem \ref{th1}}\label{sec3}
\subsection{Construction of Initial Data}\label{sec3.1}
We need to introduce smooth, radial cut-off functions to localize the frequency region. Precisely,
let $\widehat{\phi}\in \mathcal{C}^\infty_0(\mathbb{R})$ be an even, real-valued and non-negative function on $\R$ and satisfy
\begin{numcases}{\widehat{\phi}(\xi)=}
1,&if $|\xi|\leq \frac{1}{4}$,\nonumber\\
0,&if $|\xi|\geq \frac{1}{2}$.\nonumber
\end{numcases}
In \cite{Li3}, we have verified that for $f_n=\phi(x)\cos \big(\frac{17}{12}2^{n}x\big)$ and $n\geq 2$,
\begin{numcases}{\Delta_j(f_n)=}
f_n, &if $j=n$,\nonumber\\
0, &if $j\neq n$.\nonumber
\end{numcases}
We can obtain the similar result:
\begin{lemma}\label{le4}
Let $4\leq k,n\in \mathbb{N}^+$. Define the function $g^{k}_{i,n}(x)$ by
$$g^{k}_{i,n}(x):=\phi(x)\cos \bi(\frac{17}{12}\big(2^{kn}\pm2^{ki}\big)x\bi)\quad\text{with}\quad 0\leq i\leq n-1.$$
Then we have
\begin{numcases}
{\Delta_j(g^k_{i,n})=}
g^k_{i,n}, &if $j=kn$,\nonumber\\
0, &if $j\neq kn$.\nonumber
\end{numcases}
\end{lemma}
{\bf Proof.}\quad Easy computations give that
\bbal
\mathcal{F}\big(g^{k}_{i,n}\big)=2^{-1}\bi[\widehat{\phi}\bi(\xi+\frac{17}{12}\big(2^{kn}\pm2^{ki}\big)\bi)
+\widehat{\phi}\bi(\xi-\frac{17}{12}\big(2^{kn}\pm2^{ki}\big)\bi)\bi],
\end{align*}
which implies
\bbal
\mathrm{supp} \ \mathcal{F}\big(g^{k}_{i,n}\big)\subset \Big\{\xi\in\R: \ \frac{17}{12}2^{kn}-\frac{17}{12}2^{ki}-\fr12\leq |\xi|\leq \frac{17}{12}2^{kn}+\frac{17}{12}2^{ki}+\fr12\Big\},
\end{align*}
then we deduce for $k,n\geq4$
\bbal
\mathrm{supp} \ \mathcal{F}\big(g^{k}_{i,n}\big)\subset \Big\{\xi\in\R: \ \frac{33}{24}2^{kn}\leq |\xi|\leq \frac{35}{24}2^{kn}\Big\}.
\end{align*}
Notice that
\bbal
\varphi(2^{-j}\xi)\equiv 1\quad \text{for}\quad \xi\in  \mathcal{C}_j:=\Big\{\xi\in\R: \ \frac{4}{3}2^{j}\leq |\xi|\leq \frac{3}{2}2^{j}\Big\},
\end{align*}
and
\bbal
\mathcal{F}\big(\Delta_j(g^k_{i,n})\big)=\varphi(2^{-j}\cdot)\mathcal{F}\big(g^{k}_{i,n}\big),
\end{align*}
thus, for $j= kn$ we have
$$\mathcal{F}\big(\Delta_j(g^k_{i,n})\big)=\mathcal{F}\big(g^{k}_{i,n}\big).$$
Also, for $j\neq kn$ we have
$$\mathcal{F}\big(\Delta_j(g^k_{i,n})\big)=0.$$
This completes the proof of Lemma \ref{le2}.

\begin{lemma}\label{le5}
Define the initial data $u_0(x)$ as
\bbal
u_0(x):=\sum\limits^{\infty}_{n=0}2^{-kn\sigma} f^k_n(x),
\end{align*}
where
$$f^k_n(x):=\phi(x)\cos \bi(\frac{17}{12}2^{kn}x\bi),\quad n\geq 0.$$
Then for any $\sigma\in \big(2+\max\{\frac32,1+\frac1p\},+\infty\big)$ and for some $k$ large enough, we have
\bbal
&\|u_0\|_{B^{\sigma}_{p,\infty}}\leq C,\\
&\|\De_{kn}\big(u^2_0\big)\|_{L^p}\geq c2^{-kn\sigma}.
\end{align*}
\end{lemma}
{\bf Proof.}\quad By the definition of Besov space and the support of $\varphi(2^{-j}\cdot)$, using Lemma \ref{le2}, we have
\bbal
\|u_0\|_{B^{\sigma}_{p,\infty}}&=\sup_{j\geq -1}2^{\sigma j}\|\Delta_{j}u_0\|_{L^p}\\
&=\sup_{j\geq 0}\bi\|\phi(x)\cos \bi(\frac{17}{12}2^{j}x\bi)\bi\|_{L^p}\\
&\leq C.
\end{align*}
Notice that the simple fact
$$\cos(\mathbf{A}+\mathbf{B})+\cos(\mathbf{A}-\mathbf{B})=2\cos \mathbf{A} \cos \mathbf{B}$$
and
$$\sum^{\infty}_{n=0}\sum^{\infty}_{i=0,i\neq n}\mathbf{A}_n\mathbf{A}_i=2\sum^{\infty}_{n=0}\sum^{n-1}_{i=0}\mathbf{A}_n\mathbf{A}_i,$$
then direct computations give
\bbal
u^2_0(x)&=\fr12\sum^{\infty}_{n=0}2^{-2kn\sigma}\phi^2(x)
+\fr12\sum^{\infty}_{n=0}2^{-2kn\sigma}\phi^2(x)\cos\bi(\frac{17}{12}2^{kn+1}x\bi)\\
&\quad +\sum^{\infty}_{n=1}\sum^{n-1}_{i=0}2^{-k(n+i)\sigma}\phi^2(x)\Big[
\cos\bi(\frac{17}{12}(2^{kn}-2^{ik})x\bi)+
\cos\bi(\frac{17}{12}(2^{kn}+2^{ik})x\bi)\Big].
\end{align*}
Using Lemma \ref{le4} yields
\bbal
\De_{kn}\big(u^2_0\big)&=2^{-kn\sigma}\phi^2(x)\Big[
\cos\bi(\frac{17}{12}(2^{kn}-1)x\bi)+
\cos\bi(\frac{17}{12}(2^{kn}+1)x\bi)\Big]
\\&\quad +\sum^{n-1}_{i=1}2^{-k(n+i)\sigma}\phi^2(x)\Big[
\cos\bi(\frac{17}{12}(2^{kn}-2^{ik})x\bi)+
\cos\bi(\frac{17}{12}(2^{kn}+2^{ik})x\bi)\Big]\\&=:\mathrm{I}_1+\mathrm{I}_2,
\end{align*}
where we denote
\bbal
&\mathrm{I}_1:=2\cdot2^{-kn\sigma}\phi^2(x)
\cos\bi(\frac{17}{12}2^{kn}x\bi)\cos\bi(\frac{17}{12}x\bi),\\
&\mathrm{I}_2:=2\sum^{n-1}_{i=1}2^{-k(n+i)\sigma}\phi^2(x)
\cos\bi(\frac{17}{12}2^{kn}x\bi)\cos\bi(\frac{17}{12}2^{ik}x\bi).
\end{align*}
For the first term $\mathrm{I}_1$, we have
\bal\label{l0}
\|\mathrm{I}_1\|_{L^p}&\geq 2^{-kn\sigma}\bi\|\phi^2(x)\cos\bi(\frac{17}{12}x\bi) \cos\bi(\frac{17}{12}2^{kn}x\bi) \bi\|_{L^p}.
\end{align}
Since $\phi^2(x)\cos(\frac{17}{12}x)$ is a real-valued and continuous function on $\R$, then there exists some $\delta>0$ such that
$$\phi^2(x)\cos\bi(\frac{17}{12}x\bi)\geq \frac{\phi^2(0)}{2}\quad\text{ for any }  x\in B_{\delta}(0).$$
Thus, we have
\bal\label{l1}
\bi\|\phi^2(x)\bi(\frac{17}{12}x\bi)\cos \bi(\fr{17}{12}2^{kn}x\bi)\bi\|^p_{L^p}&\geq \frac{\delta}{2}\phi^2(0)\frac{1}{\lambda_n}\int^{\lambda_n}_{0}|\cos x|^p\dd x\quad\text{with}\quad \lambda_n:=\fr{17}{12}\delta2^{kn}.
\end{align}
Due to the fact
\bbal
\lim_{n\rightarrow \infty}\frac{1}{\lambda_n}\int_0^{\lambda_n}|\cos x|^p\dd x=\frac{1}{\pi}\int^\pi_0|\cos x|^p\dd x,
\end{align*}
then there exists a positive integer number $N$ such that for $n>N$
\bbal
\frac{1}{\lambda_n}\int_0^{\lambda_n}|\cos x|^p\dd x\geq\frac{1}{2\pi}\int^\pi_0|\cos x|^p\dd x,
\end{align*}
from which, \eqref{l1} reduces to
\bal\label{l2}
\bi\|\phi^2(x)\bi(\frac{17}{12}x\bi)\cos \bi(\fr{17}{12}2^{kn}x\bi)\bi\|_{L^p}&\geq c\big(p,\delta,\phi(0)\big).
\end{align}
Then we obtain from \eqref{l0}
\bal\label{l3}
\|\mathrm{I}_1\|_{L^p}&\geq c2^{-kn\sigma}.
\end{align}
For the second term $\mathrm{I}_2$, it is not hard to deduce that
\bal\label{l4}
\|\mathrm{I}_2\|_{L^p}\leq C\|\phi\|^2_{L^{2p}}\sum^{n-1}_{i=1}2^{-k(n+i)\sigma}\leq C2^{-k(n+1)\sigma}.
\end{align}
Combining \eqref{l3} and \eqref{l4} yields that
\bbal
\|\De_{kn}\big(u^2_0\big)\|_{L^p}\geq(c-C2^{-k\sigma})2^{-kn\sigma}.
\end{align*}
We choose $k\geq4$ such that $c-C2^{-k\sigma}>0$ and then finish the proof of Proposition \ref{pro3.2}.
\subsection{Error Estimates}\label{sec3.2}
\begin{proposition}\label{pro3.1}
Let $s=\sigma-2$ and $u_0\in B^{\sigma}_{p,\infty}$. Assume that $u\in L^\infty_TB^{\sigma}_{p,\infty}$ be the solution of the Cauchy problem \eqref{ch}, we have
\bbal
&\|\mathbf{S}_{t}(u_0)-u_0\|_{B^{s-1}_{p,\infty}}\leq Ct\|u_0\|_{B^{s-1}_{p,\infty}}\|u_0\|_{B^{s}_{p,\infty}},
\\&\|\mathbf{S}_{t}(u_0)-u_0\|_{B^{s}_{p,\infty}}\leq Ct\big(\|u_0\|^{2}_{B^{s}_{p,\infty}}+\|u_0\|_{B^{s-1}_{p,\infty}}\|u_0\|_{B^{s+1}_{p,\infty}}\big),
\\&\|\mathbf{S}_{t}(u_0)-u_0\|_{B^{s+1}_{p,\infty}}\leq Ct\big(\|u_0\|_{B^{s}_{p,\infty}}\|u_0\|_{B^{s+1}_{p,\infty}}+\|u_0\|_{B^{s-1}_{p,\infty}}\|u_0\|_{B^{s+2}_{p,\infty}}\big).
\end{align*}
\end{proposition}
{\bf Proof.}\quad For simplicity, we denote $u(t):=\mathbf{S}_t(u_0)$ here and in what follows. Due to the fact $B^s_{p,\infty}\hookrightarrow\rm Lip$, we know that there exists a positive time $T=T(\|u_0\|_{B^s_{p,\infty}})$ such that
\bal\label{s}
\|u(t)\|_{L^\infty_TB^s_{p,\infty}}\leq C\|u_0\|_{B^s_{p,\infty}}\leq C.
\end{align}
Moreover, for $\gamma>\fr12$, we have\bal\label{s1}
\|u(t)\|_{L^\infty_TB^\gamma_{p,\infty}}\leq C\|u_0\|_{B^\gamma_{p,\infty}}.
\end{align}
By the Mean Value Theorem, we obtain from \eqref{ch} that
\bbal
\|u(t)-u_0\|_{B^s_{p,\infty}}
&\leq \int^t_0\|\pa_\tau u\|_{B^s_{p,\infty}} \dd\tau
\\&\leq \int^t_0\|\mathbf{P}(u)\|_{B^s_{p,\infty}} \dd\tau+ \int^t_0\|u\pa_xu\|_{B^s_{p,\infty}} \dd\tau
\\&\leq Ct\big(\|u\|^2_{L_t^\infty B^{s}_{p,\infty}}+\|u\|_{L_t^\infty L^\infty}\|u_x\|_{L_t^\infty B^{s}_{p,\infty}}\big)
\\&\leq Ct\big(\|u\|^2_{L_t^\infty B^{s}_{p,\infty}}+\|u\|_{L_t^\infty B^{s-1}_{p,\infty}}\|u_x\|_{L_t^\infty B^{s}_{p,\infty}}\big)
\\&\leq Ct\big(\|u_0\|^2_{B^{s}_{p,\infty}}+\|u_0\|_{B^{s-1}_{p,\infty}}\|u_0\|_{B^{s+1}_{p,\infty}}\big),
\end{align*}
where we have used that $B_{p, \infty}^{s-1}\hookrightarrow L^\infty$ with $s-1>\max\{\frac{1}{p}, \frac{1}{2}\}$.

Following the same procedure as above, by Lemmas \ref{le1} and \ref{le2}, we have
\bbal
\|u(t)-u_0\|_{B^{s-1}_{p,\infty}}
&\leq \int^t_0\|\pa_\tau u\|_{B^{s-1}_{p,\infty}} \dd\tau
\\&\leq \int^t_0\|\mathbf{P}(u)\|_{B^{s-1}_{p,\infty}} \dd\tau+ \int^t_0\|u\pa_xu\|_{B^{s-1}_{p,\infty}} \dd\tau
\\&\leq Ct\|u\|_{L_t^\infty B^{s-1}_{p,\infty}}\|u\|_{L_t^\infty B^{s}_{p,\infty}}
\\&\leq Ct\|u_0\|_{B^{s-1}_{p,\infty}}\|u_0\|_{B^{s}_{p,\infty}},
\end{align*}
and
\bbal
\|u(t)-u_0\|_{B^{s+1}_{p,\infty}}
&\leq \int^t_0\|\pa_\tau u\|_{B^{s+1}_{p,\infty}} \dd\tau
\\&\leq \int^t_0\|\mathbf{P}(u)\|_{B^{s+1}_{p,\infty}} \dd\tau+ \int^t_0\|u\pa_xu\|_{B^{s+1}_{p,\infty}} \dd\tau
\\&\leq Ct\big(\|u\|_{L_t^\infty B^{s}_{p,\infty}}\|u\|_{L_t^\infty B^{s+1}_{p,\infty}}
+\|u\|_{L_t^\infty B^{s-1}_{p,\infty}}\|u\|_{L_t^\infty B^{s+2}_{p,\infty}}\big)
\\&\leq Ct\big(\|u_0\|_{B^{s}_{p,\infty}}\|u_0\|_{B^{s+1}_{p,\infty}}+\|u_0\|_{B^{s-1}_{p,\infty}}\|u_0\|_{B^{s+2}_{p,\infty}}\big).
\end{align*}
Thus, we finish the proof of Proposition \ref{pro3.1}.

\begin{proposition}\label{pro3.2}
Let $s=\sigma-2$ and $u_0\in B^{\sigma}_{p,\infty}$. Assume that $u\in L^\infty_TB^{\sigma}_{p,\infty}$ be the solution of the Cauchy problem \eqref{ch}, we have
\bbal
\|\mathbf{w}(t,u_0)\|_{B^{s}_{p,\infty}}\leq Ct^2\big(\|u_0\|^3_{B^s_{p,\infty}}+\|u_0\|_{B^{s-1}_{p,\infty}}\|u_0\|_{B^{s}_{p,\infty}}\|u_0\|_{B^{s+1}_{p,\infty}}
+\|u_0\|^2_{B^{s-1}_{p,\infty}}\|u_0\|_{B^{s+2}_{p,\infty}}\big),
\end{align*}
here and in what follows we denote
$$\mathbf{w}(t,u_0):=\mathbf{S}_{t}(u_0)-u_0-t\mathbf{v}_0 \quad \text{with}\quad\mathbf{v}_0:=\mathbf{P}(u_0)-u_0\pa_x u_0.$$
In particular, we have
\bbal
\|\mathbf{w}(t,u_0)\|_{B^{\sigma-2}_{p,\infty}}\leq C\big(\|u_0\|_{B^{\sigma}_{p,\infty}}\big)t^{2}.
\end{align*}
\end{proposition}
{\bf Proof.}\quad By the Mean Value Theorem and \eqref{ch}, then using Lemmas \ref{le1} and \ref{le2}, we obtain that
\bbal
\|\mathbf{w}(t,u_0)\|_{B^s_{p,\infty}}
&\leq \int^t_0\|\pa_\tau u-\mathbf{v}_0\|_{B^s_{p,\infty}} \dd\tau
\\&\leq \int^t_0\|\mathbf{P}(u)-\mathbf{P}(u_0)\|_{B^s_{p,\infty}} \dd\tau+\int^t_0\|u\pa_xu-u_0\pa_xu_0\|_{B^s_{p,\infty}} \dd\tau
\\&\leq \int^t_0\|u(\tau)-u_0\|_{B^s_{p,\infty}} \|u_0\|_{B^s_{p,\infty}}\dd\tau+\int^t_0\|u(\tau)-u_0\|_{B^{s-1}_{p,\infty}} \|u(\tau)\|_{B^{s+1}_{p,\infty}} \dd\tau
\\&\quad \ + \int^t_0\|u(\tau)-u_0\|_{B^{s+1}_{p,\infty}}  \|u_0\|_{B^{s-1}_{p,\infty}}\dd \tau
\\&\leq Ct^2\big(\|u_0\|^3_{B^s_{p,\infty}}+\|u_0\|_{B^{s-1}_{p,\infty}}\|u_0\|_{B^{s}_{p,\infty}}\|u_0\|_{B^{s+1}_{p,\infty}}
+\|u_0\|^2_{B^{s-1}_{p,\infty}}\|u_0\|_{B^{s+2}_{p,\infty}}\big),
\end{align*}
where we have used Proposition \ref{pro3.1} in the last step.

Thus, we complete the proof of Proposition \ref{pro3.2}.

Now we present the proof of Theorem \ref{th1}.\\
{\bf Proof of Theorem \ref{th1}.}\quad
Using Proposition \ref{pro3.2} and Lemma \ref{le5}, we have
\bbal
\|\mathbf{S}_{t}(u_0)-u_0\|_{B^\sigma_{p,\infty}}
&\geq2^{{kn\sigma}}\big\|\De_{kn}\big(\mathbf{S}_{t}(u_0)-u_0\big)\big\|_{L^p}=2^{{kn\sigma}}\big\|\De_{kn}\big(t\vv_0+\mathbf{w}(t,u_0)\big)\big\|_{L^p}\\
&\geq t2^{{kn\sigma}}\|\De_{kn}\big(\mathbf{v}_0\big)\|_{L^p}-2^{{2kn}}2^{{kn(\sigma-2)}}\big\|\De_{kn}\big(\mathbf{w}(t,u_0)\big)\big\|_{L^p}\\
&\geq t2^{{kn}(\sigma+1)}\|\De_{kn}\big(u^2_0\big)\|_{L^p}-
Ct\|u^2_0,(\pa_xu_0)^2\|_{B^{\sigma-1}_{p,\infty}}-C2^{2{kn}}\|\mathbf{w}(t,u_0)\|_{B^{\sigma-2}_{p,\infty}}
\\&\geq t2^{{kn}(\sigma+1)}\|\De_{kn}\big(u^2_0\big)\|_{L^p}-
Ct-C2^{2{kn}}t^2\\
&\geq ct2^{{kn}}-Ct-C2^{2{kn}}t^2.
\end{align*}
Then, for $k$, taking large $n>N$ such that $c2^{{kn}}\geq 2C$, we have
\bbal
\|\mathbf{S}_{t}(u_0)-u_0\|_{B^\sigma_{p,\infty}}\geq ct2^{{kn}}-C2^{2{kn}}t^2.
\end{align*}
Thus, picking $t2^{kn}\approx\ep$ with small $\ep$, we have
\bbal
\|\mathbf{S}_{t}(u_0)-u_0\|_{B^\sigma_{p,\infty}}\geq c\ep-C\ep^2\geq c_1\ep.
\end{align*}
This completes the proof of Theorem \ref{th1}.

\section{Proof of Theorem \ref{th2}}\label{sec4}
Firstly, we need to construct the initial data for the Novikov equation \eqref{nov}.
\begin{lemma}\label{le6}
Define the initial data $u_0(x)$ as
\bal\label{def2}
&\widehat{u_0}(\xi)=(1+|\xi|)^{-\sigma-\frac12}\quad\text{for all}\quad \xi\in\R
\end{align}
Then we have for $j\geq2$
\bbal
&\|u_0\|_{B^{\sigma}_{2,\infty}}\leq C,\\
&\|\De_j\big(u^3_0\big)\|_{L^2}\geq C2^{-\sigma j},
\end{align*}
where $C$ is some positive constant.
\end{lemma}
{\bf Proof.}\quad By the definition the Besov space and the support of $\varphi(2^{-j}\cdot)$ , we have
\bbal
\|u_0\|_{B^{\sigma}_{2,\infty}}&=\sup_{j\geq-1}2^{j\sigma}\|\Delta_{j}u_0\|_{L^2}=\sup_{j\geq-1}2^{j\sigma}\|\varphi(2^{-j}\xi)\widehat{u_0}(\xi)\|_{L^2}\leq C,
\end{align*}
which means that $u_0\in B^{\sigma}_{2,\infty}$.

In view of \eqref{def2}, by direct computations, one has
\bal\label{h1}
\mathcal{F}(u^2_0)(\xi)&=\widehat{u_0}\star\widehat{u_0}(\xi)=\int_{\R}\widehat{u_0}(\eta)\widehat{u_0}(\xi-\eta)\dd \eta\nonumber\\
&=\int_{\R}(1+|\eta|)^{-\sigma-\frac12}
\cdot(1+|\xi-\eta|)^{-\sigma-\frac12}\dd \eta\nonumber\\
&\geq\int_{|\eta|\leq 1}(1+|\eta|)^{-\sigma-\frac12}
\cdot(1+|\xi|+|\eta|)^{-\sigma-\frac12}\dd \eta\nonumber\\
&\geq (2+|\xi|)^{-\sigma-\frac12}\int_{|\eta|\leq 1}(1+|\eta|)^{-\sigma-\frac12}\dd \eta\nonumber\\
&\geq c(\sigma)\big(2+|\xi|\big)^{-\sigma-\frac12},
\end{align}
where we have used
$$\int_{|\eta|\leq 1}(1+|\eta|)^{-\sigma-\frac12}\dd \eta=2\int_{0}^{1}(1+\tau)^{-\sigma-\frac12}\dd \tau=\frac{4\big(1-2^{-\sigma+\frac12}\big)}{2\sigma-1}=:c(\sigma).$$
Furthermore, by \eqref{h1}, we obtain that
\bal\label{h2}
\mathcal{F}(u^3_0)(\xi)&=\widehat{u_0}\star\widehat{u^2_0}(\xi)=\int_{\R}\widehat{u_0}(\eta)\widehat{u^2_0}(\xi-\eta)\dd \eta\nonumber\\
&\geq c(\sigma)\int_{|\eta|\leq1}(1+|\eta|)^{-\sigma-\frac12}
\cdot(2+|\xi-\eta|)^{-\sigma-\frac12}\dd \eta\nonumber\\
&\geq c^2(\sigma)\big(3+|\xi|\big)^{-\sigma-\frac12}.
\end{align}
Then, we have for $j\geq2$
\bbal
\|\De_j\big(u^3_0\big)\|_{L^2}^2&=\int_{\fr342^j\leq|\xi|\leq\fr832^j}\varphi^2(2^{-j}\cdot)|\mathcal{F}\big(u^3_0\big)(\cdot)|^2\dd\xi\\
&\geq c^4(\sigma)\Big(3+\fr832^j\Big)^{-2\sigma-1}2^j\int_{\fr34\leq|\xi|\leq\fr83}\varphi^2\dd\xi\\
&\geq c2^{-2\sigma j}\int_{\fr43\leq|\xi|\leq\fr32}\dd\xi\\
&=c2^{-2\sigma j}.
\end{align*}
This completes the proof of Lemma \ref{le6}.

Following the similar argument as in Section \ref{sec3.2}, we can establish the following Propositions.
\begin{proposition}\label{pro4.1}
Let $s=\sigma-2$ and $u_0\in B^{\sigma}_{2,\infty}$. Assume that $u\in L^\infty_TB^{\sigma}_{2,\infty}$ be the solution of the Cauchy problem \eqref{nov}, we have
\bbal
&\|\mathbf{S}_{t}(u_0)-u_0\|_{B^{s-1}_{2,\infty}}\leq Ct\big(\|u_0\|^2_{B^{s-1}_{2,\infty}}\|u_0\|_{B^{s}_{2,\infty}}+\|u_0\|^3_{B^{s}_{2,\infty}}\big),
\\&\|\mathbf{S}_{t}(u_0)-u_0\|_{B^{s}_{2,\infty}}\leq Ct\big(\|u_0\|^{3}_{B^{s}_{2,\infty}}+\|u_0\|^2_{B^{s-1}_{2,\infty}}\|u_0\|_{B^{s+1}_{2,\infty}}\big),
\\&\|\mathbf{S}_{t}(u_0)-u_0\|_{B^{s+1}_{2,\infty}}\leq Ct\big(\|u_0\|^2_{B^{s}_{2,\infty}}\|u_0\|_{B^{s+1}_{2,\infty}}+\|u_0\|^2_{B^{s-1}_{2,\infty}}\|u_0\|_{B^{s+2}_{2,\infty}}\big).
\end{align*}
Furthermore, we have
\bbal
\|\mathbf{\widetilde{w}}(t,u_0)\|_{B^{s}_{2,\infty}}\leq Ct^2\big(\|u_0\|^3_{B^s_{2,\infty}}+\|u_0\|^2_{B^{s-1}_{2,\infty}}\|u_0\|_{B^{s+1}_{2,\infty}}
+\|u_0\|^4_{B^{s-1}_{2,\infty}}\|u_0\|_{B^{s+2}_{2,\infty}}\big),
\end{align*}
here and in what follows we denote
$$\mathbf{\widetilde{w}}(t,u_0):=\mathbf{S}_{t}(u_0)-u_0-t\mathbf{\widetilde{v}}_0\quad\text{with}\quad \mathbf{\widetilde{v}}_0:=\mathbf{Q}(u_0)-u^2_0\pa_x u_0.$$
In particular, we have
\bbal
\|\mathbf{\widetilde{w}}(t,u_0)\|_{B^{\sigma-2}_{2,\infty}}\leq C\big(\|u_0\|_{B^{\sigma}_{2,\infty}}\big)t^{2}.
\end{align*}
\end{proposition}

Now we present the proof of Theorem \ref{th2}.\\
{\bf Proof of Theorem \ref{th2}.}\quad
Using Proposition \ref{pro4.1} and Lemma \ref{le6}, we have
\bbal
\|\mathbf{S}_{t}(u_0)-u_0\|_{B^\sigma_{2,\infty}}
&\geq2^{{j\sigma}}\big\|\De_{j}\big(\mathbf{S}_{t}(u_0)-u_0\big)
\big\|_{L^2}=2^{{j\sigma}}\big\|\De_{j}
\big(t{\mathbf{\widetilde{v}}}_0+\mathbf{\widetilde{w}}(t,u_0)\big)\big\|_{L^2}\\
&\geq t2^{{j\sigma}}\|\De_{j}\big(\mathbf{\widetilde{v}}_0\big)\|_{L^2}-
2^{{2j}}2^{{j(\sigma-2)}}
\big\|\De_{j}\big(\mathbf{{\widetilde{w}}}(t,u_0)\big)\big\|_{L^2}\\
&\geq t2^{{j}(\sigma+1)}\|\De_{j}\big(u^3_0\big)\|_{L^2}-
Ct\|u^3_0,u_0(\pa_xu_0)^2,(\pa_xu_0)^3\|_{B^{\sigma-1}_{2,\infty}}\\
&~~~~-C2^{2{j}}\|\mathbf{\widetilde{w}}(t,u_0)\|_{B^{\sigma-2}_{2,\infty}}
\\&\geq t2^{{j}(\sigma+1)}\|\De_{j}\big(u^3_0\big)\|_{L^2}-
Ct-C2^{2{j}}t^2\\
&\geq ct2^{{j}}-Ct-C2^{2{j}}t^2.
\end{align*}
Then, taking large $j>N$ such that $c2^{{j}}\geq 2C$, we have
\bbal
\|\mathbf{S}_{t}(u_0)-u_0\|_{B^\sigma_{2,\infty}}\geq ct2^{{j}}-C2^{2{j}}t^2.
\end{align*}
Thus, picking $t2^{j}\approx\ep$ with small $\ep$, we have
\bbal
\|\mathbf{S}_{t}(u_0)-u_0\|_{B^\sigma_{p,\infty}}\geq c\ep-C\ep^2\geq c_2\ep.
\end{align*}
This completes the proof of Theorem \ref{th2}.

\vspace*{1em}
\noindent\textbf{Acknowledgements.} J. Li is supported by the National Natural Science Foundation of China (No.11801090) and Postdoctoral Science Foundation of China (2020T130129 and 2020M672565). Y. Yu is supported by the Natural Science Foundation of Anhui Province (No.1908085QA05) and the PhD Scientific Research Start-up Foundation of Anhui Normal University. W. Zhu is partially supported by the National Natural Science Foundation of China (No.11901092) and Natural Science Foundation of Guangdong Province (No.2017A030310634).

\vspace*{1em}
\noindent\textbf{Conflict of interest}
The authors declare that they have no conflict of interest.


\begin{thebibliography}{99}
\linespread{0}\addtolength{\itemsep}{-1.0ex}

\bibitem{B.C.D} H. Bahouri, J. Y. Chemin, R. Danchin, Fourier Analysis and Nonlinear Partial Differential Equations, Grundlehren der Mathematischen Wissenschaften, Springer, Heidelberg, 2011.

\bibitem{Byers} P. Byers, Existence time for the Camassa-Holm equation and the critical Sobolev index, Indiana Univ. Math. J. 55, 941-954 (2006).


\bibitem{Constantin} A. Constantin, Existence of permanent and breaking waves for a shallow water equation: ageometric approach, Ann. Inst. Fourier, 50 (2000), 321-362.

\bibitem{Constantin-E} A. Constantin, The Hamiltonian structure of the Camassa-Holm equation, Exposition. Math., 15 (1997), 53-85.

\bibitem{Constantin-P} A. Constantin, On the scattering problem for the Camassa-Holm equation, R. Soc. Lond. Proc. Ser. A Math. Phys. Eng. Sci., 457 (2001), 953-970.

\bibitem{Constantin-I} A. Constantin, The trajectories of particles in Stokes waves, nvent. Math., 166 (2006), 523-535.

\bibitem{Escher1} A. Constantin, J. Escher, Global existence and blow-up for a shallow water equation, Ann. Scuola Norm. Sup. Pisa Cl. Sci. (4),  26 (1998), 303-328.

\bibitem{Escher2} A. Constantin, J. Escher, Well-posedness, global existence, and blowup phenomena for a periodic quasi-linear hyperbolic equation, Comm. Pure Appl. Math., 51 (1998), 475-504.

\bibitem{Escher3} A. Constantin, J. Escher, Wave breaking for nonlinear nonlocal shallow water equations, Acta Math., 181 (1998), 229-243.

\bibitem{Escher4} A. Constantin, J. Escher, Particle trajectories in solitary water waves, Bull. Amer. Math. Soc., 44 (2007), 423-431.

\bibitem{Escher5} A. Constantin, J. Escher, Analyticity of periodic traveling free surface water waves with vorticity, Ann. Math., 173 (2011), 559-568.

\bibitem{Constantin.Strauss} A. Constantin, W. A. Strauss, Stability of peakons, Comm. Pure Appl. Math., 53 (2000), 603-610.

\bibitem{Camassa} R. Camassa, D. Holm, An integrable shallow water equation with peaked solitons, Phys. Rev. Lett., 71 (1993), 1661-1664.

\bibitem{DP} A. Degasperis, D. Holm, A. Hone, A new integral equation with peakon solutions. Theoret. Math.
Phys. 133 (2002), 1463-1474.

\bibitem{d1} R. Danchin, A few remarks on the Camassa-Holm equation, Differential Integral Equations, 14 (2001), 953-988.

\bibitem{d3} R. Danchin, A note on well-posedness for Camassa-Holm equation, J. Differential Equations, 192 (2003), 429-444.

\bibitem{Fokas} A. Fokas, B. Fuchssteiner, Symplectic structures, their B\"{a}cklund transformation and hereditary symmetries, Phys. D, 4 (1981/82), 47--66.

\bibitem{Guo-Yin} Z. Guo, X. Liu, L. Molinet, Z. Yin, Ill-posedness of the Camassa-Holm and related equations in the critical space, J. Differential Equations, 266 (2019), 1698-1707.

\bibitem{16H} A. Himonas, K. Grayshan, C. Holliman, Ill-Posedness for the b-Family of Equations, J. Nonlinear Sci. 26 (2016), 1175-1190.

\bibitem{H-H} A. Himonas, C. Holliman, The Cauchy problem for the Novikov equation, Nonlinearity, 25 (2012), 449-479.

\bibitem{HHK2018} A. Himonas, C. Holliman, C. Kenig, Construction of 2-peakon solutions and ill-posedness for the
Novikov equation, SIAM J. Math. Anal, 50(3) (2018), 2968-3006.
\bibitem{H-K} A. Himonas, C. Kenig, Non-uniform dependence on initial data for the CH equation on the line, Differential Integral Equations, 22 (2009), 201--224.

\bibitem{H-K-M} A. Himonas, C. Kenig, G. Misio{\l}ek, Non-uniform dependence for the periodic CH equation, Commun. Partial Diff. Eqns, 35 (2010), 1145--1162.

\bibitem{H-M} A. Himonas, G. Misio{\l}ek, High-frequency smooth solutions and well-posedness of the Camassa-Holm
equation, Int. Math. Res. Not., 51 (2005), 3135--3151.

\bibitem{03A} D. Holm, M. Staley, Wave structures and nonlinear balances in a family of 1+1 evolutionary PDEs. Phys.
Lett. A 308(5-6) (2003), 437-444.

\bibitem{03B} D. Holm, M. Staley, Wave structure and nonlinear balances in a family of evolutionary PDEs. SIAM J.
Appl. Dyn. Syst. 3 (2003), 323-380.

\bibitem{Li-Yin1} J. Li, Z. Yin, Remarks on the well-posedness of Camassa-Holm type equations in Besov spaces, J. Differential Equations, 261 (2016), 6125-6143.

\bibitem{Li1} J. Li, Y. Yu, W. Zhu, Non-uniform dependence on initial data for the Camassa-Holm equation in Besov spaces, J. Differential Equations, 269 (2020), 8686--8700.

\bibitem{Li2} J. Li, X. Wu, Y. Yu, W. Zhu, Non-uniform dependence on initial data for the Camassa-Holm equation in the critical Besov space. J. Math. Fluid Mech., 23:36 (2021), 11 pp.

\bibitem{Li3} J. Li, M. Li, W. Zhu, Non-uniform dependence on initial data for the Novikov equation in Besov spaces, J. Math. Fluid Mech. 22:50 (2020), 10pp.

\bibitem{li2000} Y. Li, P.J. Olver, Well-posedness and blow-up solutions for an integrable nonlinearly dispersive model wave equation, J. Differential Equations 162:1 27-63 (2000).

\bibitem{Molinet} L. Molinet, On well-posedness results for the Camassa-Holm equation on the line: a survey, J. Nonlinear Math. Phys. 11(4)(2004), 521-533.




\bibitem{GB} G. Rodr\'{i}guez-Blanco, On the Cauchy problem for the Camassa-Holm equation, Nonlinear Anal. 46 (3) (2001), 309-327.
\bibitem{Toland} J. F. Toland, Stokes waves, Topol. Methods Nonlinear Anal. 7 (1996), 1-48.

\end{thebibliography}
\end{document}